\documentclass[12pt]{article}
\usepackage{color}
\definecolor{darkblue}{rgb}{0.00,0.25,0.50}
\usepackage[colorlinks,filecolor=blue,citecolor=darkblue]{hyperref}

\usepackage{amsfonts,amssymb,amsmath,amsthm}
\usepackage{url}
\usepackage{enumerate}
\usepackage[ukrainian, russian, english]{babel}
\usepackage[cp1251]{inputenc}
\begin{document}\selectlanguage{ukrainian}
\thispagestyle{empty}

\title{}

\begin{center}
\textbf{\Large Оцінки найкращих наближень та наближень сумами Фур'є класів згорток періодичних функцій невеликої гладкості в інтегральних метриках }
\end{center}
\vskip0.5cm
\begin{center}
 Т.~А.~Степанюк\\ \emph{\small
Східноєвропейський національний університет імені Лесі
Українки, Луцьк\\}
\end{center}
\vskip0.5cm

%\address{Institute of Mathematics of the National Academy of
%Sciences of Ukraine\\ 3\\ Tereshenkivska st.\\ 01601\\ Kiev, Ukraine}

\begin{abstract}
В метриках просторів $L_{s}, \ 1<s\leq\infty$, одержано  точні за порядком
оцінки найкращих  наближень та  наближень сумами Фур'є класів згорток періодичних функцій,  що
 належать одиничній кулі
простору $L_{1}$, з твірним ядром ${\Psi_{\beta}(t)=\sum\limits_{k=1}^{\infty}\psi(k)\cos(kt-\frac{\beta\pi}{2})}$, $\beta\in \mathbb{R}$, коефіцієнти $\psi(k)$ якого такі, що добуток $\psi(n)n^{1-\frac{1}{s}}$, $1<s\leq\infty$, не може прямувати до нуля швидше за кожну степеневу функцію і, крім того, $\sum\limits_{k=1}^{\infty}\psi^{s}(k)k^{s-2}<\infty$ при $1<s<\infty$ або $\sum\limits_{k=1}^{\infty}\psi(k)<\infty$
 при $s=\infty$.

\vskip 0.5cm

In metric of spaces  $L_{s}, \ 1< s\leq\infty$, we obtain exact order estimates of best  approximations and approximations
by Fourier sums of classes of convolutions the periodic functions that  belong to unit ball of space  $L_{1}$,
with generating kernel
$\Psi_{\beta}(t)=\sum\limits_{k=1}^{\infty}\psi(k)\cos(kt-\frac{\beta\pi}{2})$, $\beta\in\mathbb{R}$, whose coefficients $\psi(k)$ are such that product  $\psi(n)n^{1-\frac{1}{s}}$, $1<s\leq\infty$,
can't tend to nought faster than every power function and besides, if  $1<s<\infty$, then $\sum\limits_{k=1}^{\infty}\psi^{s}(k)k^{s-2}<\infty$ and if $s=\infty$, then $\sum\limits_{k=1}^{\infty}\psi(k)<\infty$.
\end{abstract}

\vskip 1.5cm

%%%%%%%%%%%%%%%%%%%%%%%%%%%%%%%%%%%%%%%%%%%%%%%%%%%%%%%%%%%%%%%%%%%%%%%%%
Нехай
$L_{p}$,
$1\leq p<\infty$, --- простір $2\pi$--періодичних сумовних в $p$--му
степені на $[0,2\pi)$ функцій $f(t)$, в якому  норма задана формулою

$${\|f\|_{p}:=\Big(\int\limits_{0}^{2\pi}|f(t)|^{p}dt\Big)^{\frac{1}{p}}};$$

$L_{\infty}$ --- простір
$2\pi$--періодичних вимірних і суттєво обмежених функцій $f(t)$ з нормою
$$\|f\|_{\infty}:=\mathop{\rm{ess}\sup}\limits_{t}|f(t)|.$$

  Позначимо через
$L^{\psi}_{\beta,1}$ ---  множину функцій $f\in L_{1}$, які
 майже скрізь зображуються   за допомогою згортки
\begin{equation}\label{conv}
f(x)=\frac{a_{0}}{2}+\frac{1}{\pi}\int\limits_{-\pi}^{\pi}\Psi_{\beta}(x-t)\varphi(t)dt,
  \ \varphi\in B^{0}_{1}, \ \ \varphi\perp1, \ \ a_{0}\in\mathbb{R},
\end{equation}
де
$$
B_{1}^{0}=\left\{\varphi: \ ||\varphi||_{1}\leq 1\right\},
$$
 з сумовним ядром $\Psi_{\beta}$ вигляду
$$
\Psi_{\beta}(t)=\sum\limits_{k=1}^{\infty}\psi(k)\cos
\big(kt-\frac{\beta\pi}{2}\big), \ \psi(k)>0, \  \beta\in
    \mathbb{R}.
$$
Якщо $f$ i $\varphi$ пов'язані рівністю (\ref{conv}), то функцію $\varphi$ в цій рівності називають
$(\psi,\beta)$--похідною функції $f$ і позначають
$f_{\beta}^{\psi}$ (див., наприклад, \cite[с. 132]{Stepanets1}).

 Якщо послідовності $\psi(k)$   монотонно незростають і виконується умова
 \begin{equation}\label{cond1}
\sum\limits_{k=1}^{\infty}\psi^{s}(k)k^{s-2}<\infty, \ \ 1<s<\infty,
\end{equation}
тоді
згідно  з лемою 12.6.6 монографії \cite[с. 193]{Zigmund2}   $\Psi_{\beta}\in L_{s}$,
${1< s< \infty}$, а отже  в силу твердження 1.5.5 монографії
\cite[с. 43]{Korn} ${L^{\psi}_{\beta,1}\subset L_{s}}$, ${1< s< \infty}$.

Якщо ж $\psi(k)$
 така, що
 $$
 \sum\limits_{k=1}^{\infty}\psi(k)<\infty,
 $$
то $\Psi_{\beta}\in L_{\infty}$, і має місце включення ${L^{\psi}_{\beta,1}\subset L_{\infty}}$
(див., наприклад, твердження 1.5.5 монографії
\cite[с. 43]{Korn}).

При $\psi(k)=k^{-r}$, $r>0$ класи $L^{\psi}_{\beta,1}$ є
відомими класами  Вейля--Надя $W^{r}_{\beta,1}$, для яких при $r>1-\frac{1}{s}$  має місце включення $W^{r}_{\beta,1}\subset L_{s}$,
 $1<s\leq \infty$, $\beta\in  \mathbb{R}$.

Будемо вважати, що послідовності $\psi (k),\ k\in \mathbb{N},$ що визначають класи
 $L^{\psi}_{\beta,1}$, є звуженнями на множину
натуральних чисел  деяких додатних, неперервних, опуклих
донизу функцій $\psi(t)$ неперервного аргументу $t\geq1$, що задовольняють умову $
\lim\limits_{t\rightarrow\infty}\psi(t)=0. $  Множину всіх таких
функцій $\psi(t)$ позначатимемо через ${\mathfrak M}$.

Для класифікації функцій $\psi$ із $\mathfrak{M}$ за їх швидкістю спадання до нуля важливу роль відіграє характеристика
\begin{equation}\label{for301}
\alpha(\psi;t):=\frac{\psi(t)}{t|\psi'(t)|}, \ \ \psi'(t):=\psi'(t+0).
\end{equation}
З її допомогою з множини ${\mathfrak M}$ виділяють наступні підмножини (див., наприклад, \cite[с. 161]{Stepanets1}):
$$
\mathfrak{M}_{0}:=\left\{\psi\in \mathfrak{M}: \ \ \exists K>0 \ \ \ \  \forall t\geq1 \ \ \ \
0<K\leq \alpha(\psi;t) \right\},
$$
$$
\mathfrak{M}_{C}:=\left\{\psi\in \mathfrak{M}: \ \ \exists K_{1}, K_{2}>0 \ \ \ \ \forall t\geq1 \ \ \ \
K_{1}\leq \alpha(\psi;t)\leq K_{2}<\infty \right\},
$$
$$
\mathfrak{M}^{+}_{\infty}:=\left\{\psi\in \mathfrak{M}: \ \ \
\alpha(\psi;t)\downarrow0 \right\}.
$$

Розглянемо  величини вигляду
$$
{\cal E}_{n}(L^{\psi}_{\beta,1})_{s}=\sup\limits_{f\in
L^{\psi}_{\beta,1}}\|f(\cdot)-S_{n-1}(f;\cdot)\|_{s}, \ \ 1\leq s\leq\infty,
$$
де $S_{n-1}(f;\cdot)$ --- частинні суми Фур'є порядку $n-1$, а також найкращі  наближення класів
$L^{\psi}_{\beta,1}$ тригонометричними
поліномами  порядку не вищого за ${n-1}$, тобто величини вигляду
$$
{E}_{n}(L^{\psi}_{\beta,1})_{s}=\sup\limits_{f\in
L^{\psi}_{\beta,1}}\inf\limits_{t_{n-1}\in\mathcal{T}_{2n-1}}\|f(\cdot)-t_{n-1}(\cdot)\|_{s},
\ 1\leq s\leq \infty,
$$
де $\mathcal{T}_{2n-1}$ --- підпростір усіх тригонометричних
поліномів $t_{n-1}$ порядку не вищого за ${n-1}$.

В роботі розв'язується задача про  знаходження точних порядкових оцінок для
 величин ${\cal E}_{n}(L^{\psi}_{\beta,1})_{s}$ i ${ E}_{n}(L^{\psi}_{\beta,1})_{s}$, $1<s\leq\infty$ при певних обмеженнях на функцію $\psi$.

Для класів Вейля--Надя   $W^{r}_{\beta,1}$  порядкові оцінки величин
 ${\cal E}_{n}(W^{r}_{\beta,1})_{s}$  i
${E}_{n}(W^{r}_{\beta,1})_{s}$ при довільних $r>\frac{1}{s'}$,
$\beta\in\mathbb{R}$, ${1\leq s\leq\infty}$, $\frac{1}{s}+\frac{1}{s'}=1$ відомі (див., наприклад, \cite{T}) і мають вигляд
\begin{equation}\label{teml1}
 E_{n}(W^{r}_{\beta,1})_{s}\asymp n^{-r+\frac{1}{s'}}, \
\ \ \  1\leq s\leq\infty, \ \ \ \ r>\frac{1}{s'},
\end{equation}
$$
{\cal E}_{n}(W^{r}_{\beta,1})_{1}\asymp n^{-r}\ln n,  \ \ \ \
 \ r>0, \ \ \ \ n\in \mathbb{N}\setminus \{1\}.
$$
\begin{equation}\label{teml2}
{\cal E}_{n}(W^{r}_{\beta,1})_{s}\asymp n^{-r+\frac{1}{s'}}, \ \ \ 1<s\leq\infty, \ \ \ r>\frac{1}{s'}.
\end{equation}

В \cite{Serdyuk_grabova}
 у випадку, коли $\frac{1}{\psi(t)}$ опукла і
$\psi\in B\cap\Theta_{s'}$, $1\leq s'<\infty$, де $\Theta_{s'}$~---
множина незростаючих функцій $\psi(t)$, для яких існує
стала ${\alpha>\frac{1}{s'}}$ така, що функція $t^{\alpha}\psi(t)$
майже спадає (тобто знайдеться додатна стала $K$ для якої при будь--яких $t_{1}>t_{2}\geq 1$
${t^{\alpha}_{1}\psi(t_{1})\leq
Kt^{\alpha}_{2}\psi(t_{2})}$), а $B$ --- множина
незростаючих при $t\geq 1$ додатних функцій $\psi(t)$, для кожної з
яких можна вказати додатну сталу $K$ таку, що $
\frac{\psi(t)}{\psi(2t)}\leq K, \ \  t\geq 1 $,  показано, що існують додатні величини $K^{(1)}$, $K^{(2)}$,
   залежні лише  від $\psi$ і $s$ такі, що для довільних $1<s\leq\infty$, $\beta\in \mathbb{R}$ i $n\in \mathbb{N}$
\begin{equation}\label{gr_ser}
K^{(2)}\psi(n)n^{\frac{1}{s'}}\leq{
E}_{n}(L^{\psi}_{\beta,1})_{s}\leq{\cal
E}_{n}(L^{\psi}_{\beta,1})_{s}\leq K^{(1)}
\psi(n)n^{\frac{1}{s'}}.
\end{equation}

У випадку $\psi\in \mathfrak{M}_{\infty}^{+}$ порядкові оцінки величин ${\cal E}_{n}(L^{\psi}_{\beta,1})_{s}$
і ${ E}_{n}(L^{\psi}_{\beta,1})_{s}$, $1 \leq s\leq\infty$, були знайдені в роботах \cite{Step monog 1987}--\cite{S_S}.

В роботі
\cite{Serdyuk2013}
за умови
$\sum\limits_{k=1}^{\infty}\psi^{2}(k)<\infty$ для всіх  $\beta\in \mathbb{R}$ i $n\in\mathbb{N}$ знайдено точні значення величин
${\cal E}_{n}(L^{\psi}_{\beta,1})_{2}$, а саме встановлено рівність
$$
 {\cal E}_{n}(L^{\psi}_{\beta,1})_{2}=\frac{1}{\sqrt{\pi}}\Big(\sum\limits_{k=n}^{\infty}\psi^{2}(k)\Big)^{\frac{1}{2}}.
$$

У випадку  $\psi\in\mathfrak{M}_{C}\cup \mathfrak{M}_{\infty}^{+}$
відомі асимптотичні рівності для величин ${\cal E}_{n}(L^{\psi}_{\beta,1})_{1}$ при $n\rightarrow\infty$
(див., наприклад, \cite[с. 153]{Step monog 1987}).

При  $\psi\in {\mathfrak
M^{+}_{\infty}}$ і $\beta\in\mathbb{R}$ в \cite{serdyuk2004zbirnyk}
встановлено асимптотичні рівності і для найкращих наближень
 $E_{n}(L^{\psi}_{\beta,1})_{1}$.
Крім того, в  \cite{Serdyuk2002}--\cite{Serdyuk2005}
  отримано  точні значення величин
$E_{n}(L^{\psi}_{\beta,1})_{1}$, ${\beta\in\mathbb{R}}$, за деяких умов на послідовність $\psi(k)$.

Метою даної роботи  є  знаходження точних порядкових оцінок
 величин ${\cal E}_{n}(L^{\psi}_{\beta,1})_{s}$ i ${ E}_{n}(L^{\psi}_{\beta,1})_{s}$   у випадку, коли
 $\sum\limits_{k=1}^{\infty}\psi^{s}(k)k^{s-2}<\infty$, а функція
 \begin{equation}\label{g_s}
 g_{s'}(t):=\psi(t)t^{\frac{1}{s'}}, \ 1<s<\infty,  \ \frac{1}{s}+\frac{1}{s'}=1,
\end{equation}
 належить до множини $\mathfrak{M}_{0}$.
Крім того, в роботі знайдено точні  порядкові оцінки
 величин ${\cal E}_{n}(L^{\psi}_{\beta,1})_{\infty}$ i ${ E}_{n}(L^{\psi}_{\beta,1})_{\infty}$   у випадку, коли $\sum\limits_{k=1}^{\infty}\psi(k)<\infty$, а функція $g(t):=\psi(t)t$ належить до множини $\mathfrak{M}_{0}$.
 При цьому константи в отриманих оцінках будуть виражені через параметри  задачі в явному вигляді.

{\bf Теорема 1.} {\it
Нехай $\psi(t)t^{\frac{1}{s'}}\in\mathfrak{M}_{0}$ i
 \begin{equation}\label{condition}
 \sum\limits_{k=1}^{\infty}\psi^{s}(k)k^{s-2}<\infty, \
\end{equation}
${1<s<\infty}, \ \frac{1}{s}+\frac{1}{s'}=1$.
   Тоді   для довільних $n\in \mathbb{N}$ i $\beta\in\mathbb{R}$  мають місце співвідношення
      \begin{equation}\label{theorem_1}
 K_{\psi,s}^{(1)}\Big(\sum\limits_{k=n}^{\infty}\psi^{s}(k)k^{s-2}\Big)^{\frac{1}{s}}
 \leq{ E}_{n}(L^{\psi}_{\beta,1})_{s}\leq{\cal E}_{n}(L^{\psi}_{\beta,1})_{s} \leq
K_{\psi,s}^{(2)}\Big(\sum\limits_{k=n}^{\infty}\psi^{s}(k)k^{s-2}\Big)^{\frac{1}{s}},
\end{equation}
де $K^{(1)}_{\psi,s}$ і $K^{(2)}_{\psi,s}$ --- додатні величини
 що  залежать лише  від $\psi$ і $s$.
 }

  {\bf Доведення теореми 1.}
  Згідно з інтегральним зображенням (\ref{conv}), для довільної
функції $f\in L^{\psi}_{\beta,1}$, $\beta\in\mathbb{R}$,
майже для всіх $x\in\mathbb{R}$ справедлива  рівність
\begin{equation}\label{forrr}
f(x)-S_{n-1}(f;x)=\frac{1}{\pi}\int\limits_{-\pi}^{\pi}\Psi_{\beta,n}(x-t)f^{\psi}_{\beta}(t)dt,
\ \ \  f^{\psi}_{\beta}\in B^{0}_{1}, \ \ f^{\psi}_{\beta}\perp1,
\end{equation}
 де
 \begin{equation}\label{for6}
\Psi_{\beta,n}(t)=
\sum\limits_{k=n}^{\infty}\psi(k)\cos\big(kt-\frac{\beta\pi}{2}\big), \ \  \ n\in\mathbb{N}, \ \ \beta\in\mathbb{R}.
\end{equation}
При цьому в силу включення $g_{s'}\in\mathfrak{M}_{0}$ i умови (\ref{condition}), $\Psi_{\beta,n}\in L_{s},$ ${1<s<\infty}, \ \frac{1}{s}+\frac{1}{s'}=1$.
 Скориставшись
нерівністю (1.5.28) роботи \cite[с. 43]{Korn}, одержимо, що для довільних $1<s<\infty$
\begin{equation}\label{for34}
{\cal E}_{n}(L^{\psi}_{\beta,1})_{s}\leq
\frac{1}{\pi}\sup\limits_{f\in L_{\beta,1}^{\psi}}\|\Psi_{\beta,n}(\cdot)\|_{s}\|
f^{\psi}_{\beta}(\cdot)\|_{1}\leq\frac{1}{\pi}\|\Psi_{\beta,n}(\cdot)\|_{s}.
\end{equation}

В \cite{Serduyk_Stepaniuk} (див. формули (32) і (35)) при виконанні умови (\ref{condition}),
і умов  $g_{s'}\in\mathfrak{M}_{0}$, ${1<s<\infty}$, $\frac{1}{s}+\frac{1}{s'}=1,$ де $g_{s'}$ означається формулою (\ref{g_s}), було доведено нерівність
\begin{equation}\label{q1}
\frac{1}{\pi}\|\Psi_{\beta,n}(\cdot)\|_{s} \leq \frac{1}{\pi}\xi(s)\Big(1+ \frac{s}{\underline{\alpha}_{n}(g_{s'})}\Big)^{\frac{1}{s}}
\Big(\sum\limits_{k=n}^{\infty}\psi^{s}(k)k^{s-2}\Big)^{\frac{1}{s}},
\end{equation}
в якій
\begin{equation}\label{for100}
\xi(s):=\max\Big\{4\Big(\frac{\pi}{s-1}\Big)^{\frac{1}{s}}, \ \ 14(8\pi)^{\frac{1}{s}} s\Big\},
\end{equation}
\begin{equation}\label{k}
\underline{\alpha}_{n}(\psi):=\inf\limits_{ t\geq n}\alpha(\psi;t),
\ \ \psi\in\mathfrak{M}, \ n\in\mathbb{N},
\end{equation}
а $\alpha(\psi;t)$ означається формулою (\ref{for301}).

З нерівностей (\ref{for34}) i (\ref{q1}) отримуємо оцінку
$$
{ E}_{n}(L^{\psi}_{\beta,1})_{s}\leq{\cal E}_{n}(L^{\psi}_{\beta,1})_{s}\leq
$$
$$
\leq\frac{1}{\pi}\xi(s)\Big( \frac{\underline{\alpha}_{n}(g_{s'})+s}{\underline{\alpha}_{n}(g_{s'})}\Big)^{\frac{1}{s}}
\Big(\sum\limits_{k=n}^{\infty}\psi^{s}(k)k^{s-2}\Big)^{\frac{1}{s}}\leq
$$
\begin{equation}\label{q2}
\leq\frac{1}{\pi}\xi(s)\Big( \frac{\underline{\alpha}_{1}(g_{s'})+s}{\underline{\alpha}_{1}(g_{s'})}\Big)^{\frac{1}{s}}
\Big(\sum\limits_{k=n}^{\infty}\psi^{s}(k)k^{s-2}\Big)^{\frac{1}{s}}, \ 1<s<\infty, \frac{1}{s}+\frac{1}{s'}=1.
\end{equation}

Знайдемо  оцінку знизу для  ${ E}_{n}(L^{\psi}_{\beta,1})_{s}$, ${1<s<\infty}$. З цією метою
розглянемо згортку
\begin{equation}\label{f_m}
  f_{m}(t)=\frac{1}{\pi}\int\limits_{-\pi}^{\pi}\varphi_{m}(\tau)\Psi_{\beta}(t-\tau)d\tau,
\end{equation}
 де
\begin{equation}\label{for36}
\varphi_{m}(t):=\frac{1}{4\pi}\Big(V_{m}(t)-\frac{1}{2}\Big),
\end{equation}
а  $V_{m}(t)$ --- ядра Валле Пуссена вигляду (див. формулу (1.3.15) роботи \cite[с. 31]{Stepanets1})
\begin{equation}\label{for32}
V_{m}(t)=\frac{1}{2}+\sum\limits_{k=1}^{m}\cos kt+2\sum\limits_{k=m+1}^{2m-1}\Big(1-\frac{k}{2m}\Big)\cos
kt, \ m\in \mathbb{N}.
\end{equation}

Покажемо, що  $\|\varphi_{m}(t)\|_{1}\leq 1$, $m\in \mathbb{N}$.
Відомо, що
\begin{equation}\label{for202}
V_{m}(t)=2F_{2m-1}(t)-F_{m-1}(t),
\end{equation}
 (див., наприклад, \cite[с. 28]{T}),  де $F_{k}(t)$ --- ядра Фейєра порядку $k$
 $$
 F_{k}(t)=\frac{1}{2}+\frac{1}{k+1}\sum\limits_{\nu=0}^{k}\Big(\sum\limits_{j=1}^{\nu}\cos jt\Big), \ k\in
\mathbb{N}.
 $$
Оскільки (див., наприклад, \cite[с. 148]{Zigmund1})
\begin{equation}\label{for201}
 \|F_{k}(t)\|_{1}=\pi, \ k\in
\mathbb{N},
\end{equation}
то з (\ref{for202}) i (\ref{for201}) отримуємо
\begin{equation}\label{q3}
  \|V_{m}(t)\|_{1}\leq3\pi, \ \ m\in\mathbb{N}.
\end{equation}

Враховуючи (\ref{q3}), одержуємо
$$
\|\varphi_{m}(t)\|_{1}=\frac{1}{4\pi}\|V_{m}(t)-\frac{1}{2}\|_{1}
\leq \frac{1}{4\pi}(\|V_{m}(t)\|_{1}+\pi)\leq 1.
$$
Оскільки  $\|\varphi_{m}(t)\|_{1}\leq 1$ і $\varphi\perp1$, то  $f_{m}\in L^{\psi}_{\beta,1}, \  m\in \mathbb{N}$.  Використовуючи співвідношення  (\ref{f_m})--(\ref{for32}), а також твердження (3.7.1) з \cite[с. 134]{Stepanets1} отримаємо рівність
\begin{equation}\label{eq15}
 f_{m}(t)\!=\!\frac{1}{4\pi}\Big(\sum\limits_{k=1}^{m}\!\psi(k)\cos \Big(kt-\frac{\beta\pi}{2}\Big)+2\!\!\sum\limits_{k=m+1}^{2m-1}\!\!\!\Big(1-\frac{k}{2m}\Big)\psi(k)\cos
\Big(kt-\frac{\beta\pi}{2}\Big)\!\!\Big).
\end{equation}

Покладемо
 $$\Phi_{s}(x):=\int\limits_{x}^{\infty}\psi^{s}(t)t^{s-2}dt, \ \ 1<s<\infty,$$
 і
\begin{equation}\label{for331}
B(n)=B(\psi;s;n):=\big[\Phi_{s}^{-1}\big(\frac{1}{2n}\Phi_{s}(n)\big)\big]+1,
\end{equation}
 де $[\alpha]$ --- ціла частина дійсного числа $\alpha$, а $\Phi_{s}^{-1}$ --- функція обернена до $\Phi_{s}$.

 Розглянемо інтеграл
\begin{equation}\label{eq133}
  I_{1}=\int\limits_{-\pi}^{\pi}(f_{B(n)}(t)-t_{n-1}(t))\sum\limits_{k=n}^{\infty}
\psi^{s-1}(k)k^{s-2}\cos\Big(kt-\frac{\beta\pi}{2}\Big)dt,
\end{equation}
дe $t_{n-1}(t)\in\mathcal{T}_{2n-1}$, а функція $f_{B(n)}(t)$ означається формулою (\ref{eq15}) при $m=B(n)$.

Використавши нерівність Гельдера
(див., наприклад, \cite[с. 137]{Stepanets1}), запишемо
\begin{equation}\label{q4}
 I_{1}\leq\|f_{B(n)}(t)-t_{n-1}(t)\|_{s}\Big\|\sum\limits_{k=n}^{\infty}
\psi^{s-1}(k)k^{s-2}\cos\Big(kt-\frac{\beta\pi}{2}\Big)\Big\|_{s'}.
\end{equation}

Для оцінки норми $\Big\|\sum\limits_{k=n}^{\infty}
\psi^{s-1}(k)k^{s-2}\cos\Big(kt-\frac{\beta\pi}{2}\Big)\Big\|_{s'}$ буде корисним наступне твердження роботи \cite{Serduyk_Stepaniuk}.

{\bf Лема 1.} {\it Нехай  $1<p<\infty$
i $\{a_{k}\}_{k=1}^{\infty}$ --- монотонно незростаюча послідовність додатних чисел така, що $\sum\limits_{k=1}^{\infty}a_{k}^{p}k^{p-2}<\infty$.
Тоді для $L_{p}$--норми функції
$$
h_{\gamma,n}(x)=\sum\limits_{k=n}^{\infty}a_{k}\cos(kx+\gamma),\ \gamma\in\mathbb{R}, \ n\in\mathbb{N},
$$
має місце нерівність
$$
\|h_{\gamma,n}(x)\|_{p}\leq \xi(p)\Big(\sum\limits_{k=n}^{\infty}a_{k}^{p}k^{p-2}+a_{n}^{p}n^{p-1}\Big)^{\frac{1}{p}},
$$
де величина $\xi(p)$ означається формулою (\ref{for100}).}

Оскільки, згідно з умовою теореми,
 $g_{s'}\in\mathfrak{M}_{0}$, то функція $\psi^{s-1}(t)t^{s-2}=g_{s'}^{s-1}(t)t^{- \frac{1}{s}}$ монотонно спадає до нуля.
Тому, поклавши в умовах леми 1 $a_{k}=\psi^{s-1}(k)k^{s-2}$, ${\gamma=-\frac{\beta\pi}{2}}$,  $p=s'$,  запишемо
$$
\Big\|\sum\limits_{k=n}^{\infty}
\psi^{s-1}(k)k^{s-2}\cos\Big(kt-\frac{\beta\pi}{2}\Big)\Big\|_{s'}\leq
$$
$$
\leq \xi(s')
\bigg(\sum\limits_{k=n}^{\infty}
\big(\psi^{s-1}(k)k^{s-2}\big)^{s'}k^{s'-2}+\big(\psi^{s-1}(n)n^{s-2}\big)^{s'}n^{s'-1}\bigg)^{\frac{1}{s'}}=
$$
\begin{equation}\label{ggg5}
 =\xi(s')
\bigg(\sum\limits_{k=n}^{\infty}
\psi^{s}(k)k^{s-2}+\psi^{s}(n)n^{s-1}\bigg)^{\frac{1}{s'}}.
\end{equation}

Далі використаємо наступне твердження роботи \cite{Serduyk_Stepaniuk}.

{\bf Лема 2.} {\it
 Нехай  $\sum\limits_{k=1}^{\infty}\psi^{p'}(k)k^{p'-2}<\infty$, $1<p<\infty$, $\frac{1}{p}+\frac{1}{p'}=1$, $n\in \mathbb{N}$.
  Тоді, якщо $g_{p}\in\mathfrak{M}_{0}$, де $g_{p}(t)=\psi(t)t^{\frac{1}{p}}$, то виконується нерівність
 \begin{equation}\label{lemma_2}
\psi^{p'}(n)n^{p'-1}\leq\frac{p'}{\underline{\alpha}_{n}(g_{p})}\sum\limits_{k=n}^{\infty}\psi^{p'}(k)k^{p'-2},
\end{equation}
де величина $\underline{\alpha}_{n}(g_{p})$ означається формулою (\ref{k}).}

Застосувавши лему 2  при $p=s'$,  з (\ref{ggg5}) отримаємо
$$
\Big\|\sum\limits_{k=n}^{\infty}
\psi^{s-1}(k)k^{s-2}\cos\Big(kt-\frac{\beta\pi}{2}\Big)\Big\|_{s'}\leq
$$
\begin{equation}\label{q5}
\leq
\xi(s')\Big(1+ \frac{s}{\underline{\alpha}_{n}(g_{s'})}\Big)^{\frac{1}{s'}}
\bigg(\sum\limits_{k=n}^{\infty}
\psi^{s}(k)k^{s-2}\bigg)^{\frac{1}{s'}}.
\end{equation}

Зі співвідношень (\ref{q4}) i (\ref{q5}) отримуємо оцінку
$$
 \|f_{B(n)}(t)-t_{n-1}(t)\|_{s}\geq
$$
\begin{equation}\label{eq31}
\geq \frac{1}{\xi(s')}
 \Big(\frac{\underline{\alpha}_{n}(g_{s'})}{\underline{\alpha}_{n}(g_{s'})+s}\Big)^{\frac{1}{s'}}
\bigg(\sum\limits_{k=n}^{\infty}
\psi^{s}(k)k^{s-2}\bigg)^{-\frac{1}{s'}}
 I_{1}.
\end{equation}

Оскільки для будь--якого
 $t_{n-1}\in~ \mathcal{T}_{2n-1}$
\begin{equation}\label{jj}
\int\limits_{-\pi}^{\pi}t_{n-1}(t)\sum\limits_{k=n}^{\infty}
\psi^{s-1}(k)k^{s-2}\cos\Big(kt-\frac{\beta\pi}{2}\Big)dt=0,
\end{equation}
то в силу (\ref{eq133})
\begin{equation}\label{q7}
  I_{1}=\int\limits_{-\pi}^{\pi}f_{B(n)}(t)\sum\limits_{k=n}^{\infty}
\psi^{s-1}(k)k^{s-2}\cos\Big(kt-\frac{\beta\pi}{2}\Big)dt.
\end{equation}

Очевидно, що
\begin{equation}\label{int_riv}
  \int\limits_{-\pi}^{\pi}\cos(kt+\theta)\cos(mt+\theta)dt=
{\left\{\begin{array}{cc}
0, \ & k\neq m, \\
\pi, &
k=m, \
  \end{array} \right.}  \ \ k,m\in\mathbb{N}, \ \ \theta\in\mathbb{R.}
\end{equation}

 Використовуючи (\ref{int_riv}) при $\theta=-\frac{\beta\pi}{2}$ i (\ref{eq15}) при $m=B(n)$,
з (\ref{q7}) одержуємо
$$
I_{1}
=\frac{1}{4\pi}\int\limits_{-\pi}^{\pi}\Big( \sum\limits_{k=1}^{B(n)}\psi(k)\cos \Big(kt-\frac{\beta\pi}{2}\Big)+
$$
$$
+2\sum\limits_{k=B(n)+1}^{2B(n)-1}\Big(1-\frac{k}{2B(n)}\Big)\psi(k)\cos
\Big(kt-\frac{\beta\pi}{2}\Big) \Big)\times
$$
$$
\times\sum\limits_{k=n}^{\infty}
\psi^{s-1}(k)k^{s-2}\cos\Big(kt-\frac{\beta\pi}{2}\Big)dt=
$$
$$
=\frac{1}{4}\Big(\sum\limits_{k=n}^{B(n)}
\psi^{s}(k)k^{s-2}+2\sum\limits_{k=B(n)+1}^{2B(n)-1}\Big(1-\frac{k}{2B(n)}\Big)\psi^{s}(k)k^{s-2}\Big)>
$$
\begin{equation}\label{eq32}
  >\frac{1}{4}\sum\limits_{k=n}^{B(n)}
\psi^{s}(k)k^{s-2}=\frac{1}{4}\Big(\sum\limits_{k=n}^{\infty}
\psi^{s}(k)k^{s-2} -\sum\limits_{k=B(n)+1}^{\infty}
\psi^{s}(k)k^{s-2} \Big).
\end{equation}
 Для оцінки знизу інтеграла $I_{1}$ залишилось оцінити зверху суму $ \sum\limits_{k=B(n)+1}^{\infty}
\psi^{s}(k)k^{s-2}$.  З  (\ref{for331}) випливає
$$
 \sum\limits_{k=B(n)+1}^{\infty}\psi^{s}(k)k^{s-2}\leq\int\limits_{B(n)}^{\infty}\psi^{s}(t)t^{s-2}dt=
\Phi_{s}(B(n))<
$$
 \begin{equation}\label{eq33}
<\frac{1}{2n}\Phi_{s}(n)\leq\frac{1}{2n}\sum\limits_{k=n}^{\infty}\psi^{s}(k)k^{s-2}.
 \end{equation}

 З нерівностей (\ref{eq31}), (\ref{eq32}) і (\ref{eq33}) для довільного $t_{n-1}(t)\in\mathcal{T}_{2n-1}$ отримаємо оцінку
 $$
 \|f_{B(n)}(t)-t_{n-1}(t)\|_{s}\geq
 $$
 $$
 \geq
   \frac{1}{\xi(s')} \Big(\frac{\underline{\alpha}_{n}(g_{s'})}{\underline{\alpha}_{n}(g_{s'})+s}\Big)^{\frac{1}{s'}}
\bigg(\!\sum\limits_{k=n}^{\infty}
\psi^{s}(k)k^{s-2}\bigg)^{-\frac{1}{s'}}\!
   \frac{1}{4}\Big(1-\frac{1}{2n}\Big)\!\sum\limits_{k=n}^{\infty}
\psi^{s}(k)k^{s-2}\geq
 $$
 $$
 \geq\frac{1}{8\xi(s')} \Big(\frac{\underline{\alpha}_{n}(g_{s'})}{\underline{\alpha}_{n}(g_{s'})+s}\Big)^{\frac{1}{s'}}\bigg(\sum\limits_{k=n}^{\infty}\psi^{s}(k)k^{s-2} \bigg)^{\frac{1}{s}}\geq
 $$
  \begin{equation}\label{eq34}
\geq\frac{1}{8\xi(s')} \Big(\frac{\underline{\alpha}_{1}(g_{s'})}{\underline{\alpha}_{1}(g_{s'})+s}\Big)^{\frac{1}{s'}}\bigg(\sum\limits_{k=n}^{\infty}\psi^{s}(k)k^{s-2} \bigg)^{\frac{1}{s}}.
\end{equation}
Об'єднуючи (\ref{q2}) і (\ref{eq34}) отримуємо співвідношення (\ref{theorem_1}). Теорему 1 доведено.

Прикладами функцій $\psi$, які задовольняють умови теореми 1  є функції:
$$
 1)  \psi(t)=t^{-r}, \ {r>\frac{1}{s'}}; \ \ \ \ \ \ \ \ \ \ \ \ \ \ \ \ \ \ \ \ \ \ \ \ \ \ \ \ \ \ \ \ \ \ \ \ \ \ \ \ \ \\ \ \ \ \ \ \ \ \ \ \ \ \ \ \ \ \ \ \ \ \ \ \ \ \
$$
$$
 2) \psi(t)={t^{-\frac{1}{s'}}\ln^{-\gamma}(t+K)}, {\gamma>\frac{1}{s}}, K>0; \ \ \ \ \ \ \ \ \ \ \ \ \ \ \ \ \ \ \ \ \\ \ \ \ \ \ \ \ \ \ \ \ \ \ \ \ \ \ \ \ \
$$
\begin{equation}\label{func_3}
 3)  \psi(t)\!=\!{t^{-\frac{1}{s'}}(\ln (t\!+\!K_{1}))^{-\frac{1}{s}}
(\ln\ln (t\!+\!K_{2}))^{-\gamma}},   \gamma>\frac{1}{s},  K_{2}\geq e-1,  K_{1}>0,
\end{equation}
\noindent та інші.

{\bf Зауваження. }{\it В ході доведення теореми 1 за виконання її умов було показано, що для  довільного $n\in\mathbb{N}$ виконується більш точна, ніж (\ref{theorem_1}), оцінка:
$$
\frac{1}{8\xi(s')} \Big(\frac{\underline{\alpha}_{n}(g_{s'})}{\underline{\alpha}_{n}(g_{s'})+s}\Big)^{\frac{1}{s'}} \Big(\sum\limits_{k=n}^{\infty}\psi^{s}(k)k^{s-2}\Big)^{\frac{1}{s}}
 \leq{ E}_{n}(L^{\psi}_{\beta,1})_{s}\leq{\cal E}_{n}(L^{\psi}_{\beta,1})_{s}
\leq
$$
\begin{equation}\label{remark}
 \leq
\frac{1}{\pi}\xi(s)\Big(\frac{\underline{\alpha}_{n}(g_{s'})+s}{\underline{\alpha}_{n}(g_{s'})}\Big)^{\frac{1}{s}}
\Big(\sum\limits_{k=n}^{\infty}\psi^{s}(k)k^{s-2}\Big)^{\frac{1}{s}},
\end{equation}
де $1<s<\infty,$ $\frac{1}{s}+\frac{1}{s'}=1$, a $\xi(s')$ і $\underline{\alpha}_{n}(g_{s'})$ --- додатні величини, що означаються за
допомогою формул (\ref{for100}) і (\ref{k}) відповідно.  }

З нерівностей (\ref{remark}) випливає наступне твердження.

{\bf Наслідок 1. }{\it Нехай $r>\frac{1}{s'}$, $1<s<\infty$, $\frac{1}{s}+\frac{1}{s'}=1$. Тоді для   довільних $\beta\in\mathbb{R}$ і $n\in\mathbb{N}$ виконується оцінка: }
$$
\frac{1}{8\xi(s')} \Big(\frac{s'}{s'+s(rs'-1)}\Big)^{\frac{1}{s'}} \Big(\sum\limits_{k=n}^{\infty}  \frac{1}{k^{s(r-1)+2}}\Big)^{\frac{1}{s}}
 \leq{ E}_{n}(W^{r}_{\beta,1})_{s}\leq
$$
\begin{equation}\label{conseq1}
\leq{\cal E}_{n}(W^{r}_{\beta,1})_{s}
\leq
\frac{1}{\pi}\xi(s)\Big(\frac{s'+s(rs'-1)}{s'}\Big)^{\frac{1}{s}}
\Big(\sum\limits_{k=n}^{\infty}\frac{1}{k^{s(r-1)+2}}\Big)^{\frac{1}{s}},
\end{equation}
{\it де $\xi(s)$ --- додатня величина, що означається за
допомогою формули (\ref{for100}).}

Оскільки
$$
\frac{1}{\alpha-1}n^{1-\alpha}\leq
\sum\limits_{k=n}^{\infty}\frac{1}{k^{\alpha}}\leq\Big(1+\frac{1}{\alpha-1}\Big)n^{1-\alpha}, \ \alpha>1,
$$
то  з (\ref{conseq1}) одержуємо співвідношення
$$
\frac{1}{8\xi(s')} \Big(\frac{s'}{s'+s(rs'-1)}\Big)^{\frac{1}{s'}}
\Big(\frac{1}{s(r-1)+1}\Big)^{\frac{1}{s}}n^{-r+\frac{1}{s'}}
 \leq{ E}_{n}(W^{r}_{\beta,1})_{s}\leq
$$
\begin{equation}\label{nerW}
  \leq{\cal E}_{n}(W^{r}_{\beta,1})_{s}
\leq
\frac{1}{\pi}\xi(s)\Big(\frac{s'+s(rs'-1)}{s'}\Big)^{\frac{1}{s}}
\Big(\frac{s(r-1)+2}{s(r-1)+1}\Big)^{\frac{1}{s}}n^{-r+\frac{1}{s'}}.
\end{equation}
Нерівності (\ref{nerW}) уточнюють порядкові оцінки (\ref{teml1}) i (\ref{teml2}).

В \cite{Serduyk_Stepaniuk}  було показано, що при виконанні умов $\sum\limits_{k=1}^{\infty}\psi^{s}(k)k^{s-2}<\infty$,
$\psi(t)t^{\frac{1}{s'}}\in\mathfrak{M}_{C}$, $1<s<\infty$, $\frac{1}{s}+\frac{1}{s'}=1$  має місце порядкова оцінка
\begin{equation}\label{order}
  \sum\limits_{k=n}^{\infty}\psi^{s}(k)k^{s-2}\asymp\psi^{s}(n)n^{s-1},
\end{equation}
де під записом  ${A(n)\asymp B(n)}$, як зазвичай прийнято, будемо розуміти,
 що для додатних послідовностей $A(n)$ i $B(n)$ існують  сталі $K_{1}>0$ і
$K_{2}>0$ такі, що ${K_{1}B(n)\leq A(n)\leq K_{2}B(n)}$,
$n\in\mathbb{N}$.

Отже, з (\ref{theorem_1}) i (\ref{order}) випливає наступне твердження.

{\bf Наслідок 2. }{\it  Нехай $\sum\limits_{k=1}^{\infty}\psi^{s}(k)k^{s-2}<\infty$, ${1<s<\infty}$, $\frac{1}{s}+\frac{1}{s'}=1$, $\beta\in\mathbb{R}$ i $n\in \mathbb{N}$.
  Тоді, якщо  для функції $g_{s'}$ вигляду (\ref{g_s}) виконується включення $g_{s'}\in\mathfrak{M}_{0}$, то
 \begin{equation}\label{c1}
{ E}_{n}(L^{\psi}_{\beta,1})_{s}\asymp{\cal E}_{n}(L^{\psi}_{\beta,1})_{s}\asymp\Big(\sum\limits_{k=n}^{\infty}\psi^{s}(k)k^{s-2}\Big)^{\frac{1}{s}},
\end{equation}
а якщо ж $g_{s'}\in\mathfrak{M}_{C}$, то}
\begin{equation}\label{cc1}
{ E}_{n}(L^{\psi}_{\beta,1})_{s}\asymp{\cal E}_{n}(L^{\psi}_{\beta,1})_{s}\asymp\psi(n)n^{\frac{1}{s'}}.
\end{equation}

Порядкові оцінки (\ref{cc1})  встановлені раніше в роботі \cite{Serdyuk_grabova}.

Зауважимо, що у випадку, коли
\begin{equation}\label{lim}
g_{s'}\in\mathfrak{M}_{0}, \ \ \lim\limits_{t\rightarrow\infty}\alpha(g_{s'};t)=\infty
\end{equation}
виконується оцінка
$$
\psi(n)n^{\frac{1}{s'}}=o\bigg(\Big(\sum\limits_{k=n}^{\infty}\psi^{s}(k)k^{s-2}\Big)^{\frac{1}{s}}\bigg), \ n\rightarrow\infty,
$$
 тобто порядкові рівності (\ref{cc1}) місця не мають. Зокрема умова (\ref{lim}) виконується для функцій $\psi(t)$ вигляду (\ref{func_3}).
Наведемо наслідок з теореми 1 для згаданих функцій $\psi$.

{\bf Наслідок 3. }{\it Нехай $\psi(t)=t^{-\frac{1}{s'}}\ln^{-\frac{1}{s}} (t+K_{1})
(\ln\ln (t+K_{2}))^{-\gamma}$, ${\gamma>\frac{1}{s}}$,  $K_{2}\geq e-1$, $K_{1}>0$, $1<s<\infty$, $\frac{1}{s}+\frac{1}{s'}=1$,  $\beta\in\mathbb{R}$ i $n\geq3, \ n\in\mathbb{N}$. Тоді
$$
{E}_{n}(L^{\psi}_{\beta,1})_{s}\asymp{\cal E}_{n}(L^{\psi}_{\beta,1})_{s}\asymp (\ln\ln n)^{\frac{1}{s}-\gamma}
\asymp\psi(n)n^{\frac{1}{s'}}(\ln n )^{\frac{1}{s}}(\ln\ln n)^{\frac{1}{s}}.
$$
}
{\bf Доведення наслідку 3. } Як зазначалось вище, функції $\psi(t)$ вигляду (\ref{func_3})
задовольняють умови теореми 1. Тому, враховуючи співвідношення (42) роботи \cite{Serduyk_Stepaniuk},
неважко переконатись в справедливості співвідношення
$$
\int\limits_{n}^{\infty}\psi^{s}(t)t^{s-2}dt
\leq\sum\limits_{k=n}^{\infty}\psi^{s}(k)k^{s-2}\leq
\psi^{s}(n)n^{s-2}+\int\limits_{n}^{\infty}\psi^{s}(t)t^{s-2}dt\leq
$$
\begin{equation}\label{for110}
\leq\Big( \frac{s}{\underline{\alpha}_{n}(g_{s'})}\cdot\frac{1}{n}+1\Big)
\int\limits_{n}^{\infty}\psi^{s}(t)t^{s-2}dt.
\end{equation}
З (\ref{theorem_1}), (\ref{for110}) і того, що $\underline{\alpha}_{n}(g_{s'})>K>0$ випливає, що при  $n\geq3$
$$
{ E}_{n}(L^{\psi}_{\beta,1})_{s}
\asymp{\cal E}_{n}(L^{\psi}_{\beta,1})_{s}\asymp
\Big(\sum\limits_{k=n}^{\infty}\psi^{s}(k)k^{s-2}\Big)^{\frac{1}{s}}\asymp
$$
$$
\asymp\Big(\int\limits_{n}^{\infty}\psi^{s}(t)t^{s-2}dt\Big)^{\frac{1}{s}}=
\Big(\int\limits_{n}^{\infty}\frac{dt}{t\ln(t+K_{1})(\ln\ln(t+K_{2}))^{\gamma
s}}\Big)^{\frac{1}{s}}
\asymp
$$
$$
\asymp\Big(\int\limits_{n}^{\infty}\frac{dt}{t\ln t(\ln\ln t)^{\gamma
s}}\Big)^{\frac{1}{s}}\asymp
 (\ln\ln n)^{\frac{1}{s}-\gamma}\asymp
$$
$$
\asymp\psi(n)n^{\frac{1}{s'}}(\ln n)^{\frac{1}{s}} (\ln\ln n)^{\frac{1}{s}}, \ n\geq3.
$$
Наслідок  3 доведено.

{\bf Теорема 2.} {\it Нехай $\psi\in\mathfrak{M}$, $\sum\limits_{k=1}^{\infty}\psi(k)<\infty$, $\beta\in\mathbb{R}$ i
$\cos\frac{\beta\pi}{2}\neq0$.
  Тоді   для довільних $n\in \mathbb{N}$  мають місце нерівності
    \begin{equation}\label{theorem_2}
\frac{1}{4\pi}\Big|\cos\frac{\beta\pi}{2}\Big|\sum\limits_{k=n}^{\infty}\psi(k)\leq{\cal E}_{n}(L^{\psi}_{\beta,1})_{\infty}\leq\frac{1}{\pi}\sum\limits_{k=n}^{\infty}\psi(k).
\end{equation}
 }
 {\bf Доведення теореми 2.} Знайдемо оцінку зверху величини  ${\cal E}_{n}(L^{\psi}_{\beta,1})_{\infty}$. З урахуванням  формули (\ref{forrr}) i нерівності (1.5.28) роботи
 \cite[с. 43]{Korn}
  \begin{equation}\label{for230}
 {\cal E}_{n}(L^{\psi}_{\beta,1})_{\infty}\leq
\frac{1}{\pi}\big\|\Psi_{\beta,n}(\cdot)\big\|_{\infty}.
\end{equation}

  Оскільки
  $$
\|\Psi_{\beta,n}(t)\|_{\infty}=\Big\|\sum\limits_{k=n}^{\infty}\psi(k)\cos\Big(kt-\frac{\beta\pi}{2}\Big)\Big\|_{\infty}\leq
\sum\limits_{k=n}^{\infty}\psi(k),
 $$
 то з (\ref{for230}) одержуємо
 \begin{equation}\label{for24}
 {\cal E}_{n}(L^{\psi}_{\beta,1})_{\infty}\leq
\frac{1}{\pi}\sum\limits_{k=n}^{\infty}\psi(k), \  \beta\in\mathbb{R}, \  n\in \mathbb{N}.
\end{equation}

Знайдемо оцінку знизу  величини ${\cal E}_{n}(L^{\psi}_{\beta,1})_{\infty}$.
 Покладемо
 $$
 \Psi(x):=\int\limits_{x}^{\infty}\psi(t)dt
 $$
  і
\begin{equation}\label{for31}
D(l;n)=D(\psi;l;n):=\big[\Psi^{-1}\big(\frac{1}{2l}\Psi(n)\big)\big]+2n, \ \ l,n\in\mathbb{N}.
\end{equation}

Розглянемо функцію $f_{D(l;n)}(t)$, що означається формулою (\ref{eq15}) при $m=D(l;n)$, тобто
$$
 f_{D(l;n)}(t)=\frac{1}{4\pi}\Big(\sum\limits_{k=1}^{D(l;n)}\psi(k)\cos \Big(kt-\frac{\beta\pi}{2}\Big)+
$$
\begin{equation}\label{eeeq15}
+2\sum\limits_{k=D(l;n)+1}^{2D(l;n)-1}\Big(1-\frac{k}{2D(l;n)}\Big)\psi(k)\cos
\Big(kt-\frac{\beta\pi}{2}\Big)\Big),
\end{equation}
Як було показано
при доведенні теореми 1,
$ f_{D(l;n)}\in L^{\psi}_{\beta,1}$.
Беручи до уваги рівність (\ref{eeeq15}) та враховуючи умови теореми 2, отримуємо, що при довільних $l\in\mathbb{N}$
$$
{\cal E}_{n}(L^{\psi}_{\beta,1})_{\infty}\geq|f_{D(l;n)}(0)-S_{n-1}(f_{D(l;n)};0)|=
$$
$$
=\frac{1}{4\pi}\Big|\cos\frac{\beta\pi}{2}\Big|\Big(\sum\limits_{k=n}^{D(l;n)}\psi(k)+2\sum\limits_{k=D(l;n)+1}^{2D(l;n)-1}\Big(1-\frac{k}{2D(l;n)}\Big)\psi(k)\Big)
>
$$
\begin{equation}\label{for33}
>\frac{1}{4\pi}\Big|\cos\frac{\beta\pi}{2}\Big|\sum\limits_{k=n}^{D(l;n)}\psi(k)=
\frac{1}{4\pi}\Big|\cos\frac{\beta\pi}{2}\Big|\Big(\sum\limits_{k=n}^{\infty}\psi(k)-\sum\limits_{k=D(l;n)+1}^{\infty}\psi(k)\Big).
\end{equation}
 З (\ref{for31}) випливає, що для довільних $ l\in \mathbb{N}$
 \begin{equation}\label{for401}
 \sum\limits_{k=D(l;n)+1}^{\infty}\psi(k)\leq\int\limits_{D(l;n)}^{\infty}\psi(t)dt=\Psi(D(l;n))<
 \frac{1}{2l}\Psi(n)\leq\frac{1}{2l}\sum\limits_{k=n}^{\infty}\psi(k).
 \end{equation}
 Тому
 \begin{equation}\label{nn}
 {\cal E}_{n}(L^{\psi}_{\beta,1})_{\infty}>
 \frac{1}{4\pi}\Big|\cos\frac{\beta\pi}{2}\Big|\Big(1-\frac{1}{2l}\Big)\sum\limits_{k=n}^{\infty}\psi(k).
  \ l,n\in \mathbb{N}
 \end{equation}

Перейшовши до границі в нерівності  (\ref{nn}) при $l\rightarrow\infty$, одержимо
\begin{equation}\label{for103}
{\cal E}_{n}(L^{\psi}_{\beta,1})_{\infty}>
 \frac{1}{4\pi}\Big|\cos\frac{\beta\pi}{2}\Big|\sum\limits_{k=n}^{\infty}\psi(k).
\end{equation}
На підставі (\ref{for24}) і (\ref{for103}) отримуємо  оцінку (\ref{theorem_2}). Теорему 2 доведено.

{\bf Теорема 3.} {\it  Нехай $\sum\limits_{k=1}^{\infty}\psi(k)<\infty$, $\psi(t)=g(t)t^{-1}$, $g\in\mathfrak{M}_{0}$ і
\begin{equation}\label{cond}
 \underline{\alpha}_{1}(g)=\inf\limits_{t\geq 1}\alpha(g;t)>1.
\end{equation}

   Тоді   для довільних $\beta\in\mathbb{R}$,
$\cos\frac{\beta\pi}{2}\neq0$ i $n\in \mathbb{N}$
 \begin{equation}\label{theorem_3}
\frac{1}{48\pi}\Big|\cos\frac{\beta\pi}{2}\Big|
\Big(1-\frac{1}{\underline{\alpha}_{1}(g)}\Big)\sum\limits_{k=n}^{\infty}\psi(k)\leq{ E}_{n}(L^{\psi}_{\beta,1})_{\infty}\leq{\cal E}_{n}(L^{\psi}_{\beta,1})_{\infty}\leq\frac{1}{\pi}\sum\limits_{k=n}^{\infty}\psi(k).
\end{equation}
 }
 {\bf Доведення теореми 3.}
 В силу теореми 2 при умові $\sum\limits_{k=1}^{\infty}\psi(k)<\infty$ справедлива оцінка
\begin{equation}\label{j3}
  {E}_{n}(L^{\psi}_{\beta,1})_{\infty} \leq
  \mathcal{E}_{n}(L^{\psi}_{\beta,1})_{\infty}\leq\frac{1}{\pi}\sum\limits_{k=n}^{\infty}\psi(k), \  \beta\in\mathbb{R}, \  n\in \mathbb{N}.
\end{equation}

Знайдемо  оцінку знизу величини  ${E}_{n}(L^{\psi}_{\beta,1})_{\infty}$.
Розглянемо інтеграл
\begin{equation}\label{eq13}
  I_{2}=\int\limits_{-\pi}^{\pi}(f_{D(l;n)}(t)-t_{n-1}(t))(V_{D(l;n)}(t)-V_{n-1}(t))dt,
\end{equation}
де $t_{n-1}(t)\in\mathcal{T}_{2n-1}$, а функція $f_{D(l;n)}(t)$ та величина  $D(l;n)$  означаються формулами (\ref{eeeq15}) і (\ref{for31}) відповідно.

Використавши формулу (\ref{for32}), запишемо
$$
 V_{D(l;n)}(t)-V_{n-1}(t)=
$$
\begin{equation}\label{w1}
=\!\sum\limits_{k=n}^{D(l;n)}\!\!\cos kt+2\sum\limits_{k=D(l;n)+1}^{2D(l;n)-1}\!\!\Big(1-\frac{k}{2D(l;n)}\Big)\!\cos
kt-2\sum\limits_{k=n}^{2n-3}\!\!\Big(1-\frac{k}{2n-2}\Big)\!\cos
kt.
\end{equation}
З (\ref{w1}) випливає, що для довільного $t_{n-1}(t)\in\mathcal{T}_{2n-1}$
\begin{equation}\label{for402}
 \int\limits_{-\pi}^{\pi}t_{n-1}(t)\big( V_{D(l;n)}(t)-V_{n-1}(t)\big)dt=0  \ \ \forall  n\in\mathbb{N}.
\end{equation}

Оскільки
\begin{equation}\label{w3}
   \int\limits_{-\pi}^{\pi}\cos kt\cos\Big(mt-\frac{\beta\pi}{2}\Big)dt=
{\left\{\begin{array}{cc}
0, \ & k\neq m, \\
\pi\cos\frac{\beta\pi}{2}, &
k=m, \
  \end{array} \right.}  \ \ k,m\in\mathbb{N}, \ \beta\in\mathbb{R},
\end{equation}
то беручи до уваги (\ref{eeeq15}), (\ref{w1}) i (\ref{for402}), одержуємо
$$
 | I_{2}|=\bigg|\int\limits_{-\pi}^{\pi}f_{D(l;n)}(t)(V_{D(l;n)}(t)-V_{n-1}(t))dt\bigg|=
$$
$$
=\frac{1}{4\pi}\bigg|\int\limits_{-\pi}^{\pi}
\bigg(\sum\limits_{k=1}^{D(l;n)}\psi(k)\cos \Big(kt-\frac{\beta\pi}{2}\Big)+
2\sum\limits_{k=D(l;n)+1}^{2D(l;n)-1}\Big(1-\frac{k}{2D(l;n)}\Big)\psi(k)\cos
\Big(kt-\frac{\beta\pi}{2}\Big)\bigg)\times
$$
$$
\times\bigg(\sum\limits_{k=n}^{D(l;n)}\cos kt+2\sum\limits_{k=D(l;n)+1}^{2D(l;n)-1}\Big(1-\frac{k}{2D(l;n)}\Big)\cos
kt-
2\sum\limits_{k=n}^{2n-3}\Big(1-\frac{k}{2n-2}\Big)\cos
kt\bigg)dt\bigg|=
$$
$$
=\frac{1}{4}\Big|\cos\frac{\beta\pi}{2}\Big|
\bigg(\sum\limits_{k=n}^{D(l;n)}\psi(k)-2\sum\limits_{k=n}^{2n-3}\Big(1-\frac{k}{2n-2}\Big)\psi(k)+
4\sum\limits_{k=D(l;n)+1}^{2D(l;n)-1}\Big(1-\frac{k}{2D(l;n)}\Big)^{2}\psi(k)\bigg)>
$$
$$
 >\frac{1}{4}\Big|\cos\frac{\beta\pi}{2}\Big|
\bigg(\sum\limits_{k=n}^{D(l;n)}\psi(k)-2\sum\limits_{k=n}^{2n-3}\Big(1-\frac{k}{2n-2}\Big)\psi(k)\bigg)=
$$
\begin{equation}\label{eq16}
=\frac{1}{4}\Big|\cos\frac{\beta\pi}{2}\Big|
\bigg(\sum\limits_{k=n}^{\infty}\psi(k)-\sum\limits_{k=D(l;n)+1}^{\infty}\psi(k)-2\sum\limits_{k=n}^{2n-3}\Big(1-\frac{k}{2n-2}\Big)\psi(k)\bigg).
\end{equation}
Враховуючи монотонність функції $\psi(t)$, отримуємо
$$
 2\sum\limits_{k=n}^{2n-3}\Big(1-\frac{k}{2n-2}\Big)\psi(k)\leq\psi(n)\frac{1}{n-1}\sum\limits_{k=n}^{2n-3}(2n-2-k)=
$$
\begin{equation}\label{eq18}
=
\frac{1}{2}\psi(n)(n-2)<\frac{1}{2}\psi(n)n.
\end{equation}

Для будь--якої функції $\psi\in\mathfrak{M}$ через  $\overline{\alpha}_{n}(\psi)$, $n\in \mathbb{N}$, позначимо величину
\begin{equation}\label{kk}
\overline{\alpha}_{n}(\psi):=\sup\limits_{ t\geq n}\alpha(\psi;t),
\end{equation}
де характеристика $\alpha(\psi;t)$ означається формулою (\ref{for301}).
В прийнятих позначеннях має місце наступне твердження.

{\bf Лема 3.} {\it Нехай $\psi(t)=g(t)t^{-1}$, $g\in\mathfrak{M}_{0}$ i $\sum\limits_{k=1}^{\infty}\psi(k)<\infty$.
  Тоді   для довільних $n\in \mathbb{N}$
 \begin{equation}\label{lemma_3}
\psi(n)n\leq\frac{1}{\underline{\alpha}_{n}(g)}\sum\limits_{k=n}^{\infty}\psi(k).
\end{equation}

Якщо ж крім того $g\in\mathfrak{M}_{C}$, то}
\begin{equation}\label{lemma_3_1}
\frac{1}{\overline{\alpha}_{n}(g)}\cdot\frac{n\underline{\alpha}_{n}(g)}{1+n\underline{\alpha}_{n}(g)}\sum\limits_{k=n}^{\infty}\psi(k)
\leq\psi(n)n\leq\frac{1}{\underline{\alpha}_{n}(g)}\sum\limits_{k=n}^{\infty}\psi(k).
\end{equation}

{\bf Доведення леми 3.} Нехай $g\in\mathfrak{M}_{0}$. Покажемо справедливість нерівності (\ref{lemma_3}).
Очевидно,
\begin{equation}\label{eqq2}
\sum\limits_{k=n}^{\infty}\psi(k) \geq\int\limits_{n}^{\infty}\psi(t)dt.
\end{equation}

Оскільки функція $\psi(t)t$ монотонно спадає до нуля, то проінтегрувавши частинами інтеграл $\int\limits_{n}^{\infty}\psi(t)dt$, отримаємо
\begin{equation}\label{eqq1}
 \int\limits_{n}^{\infty}\psi(t)dt=
-\psi(n)n-\int\limits_{n}^{\infty}\psi'(t)tdt=
-\psi(n)n+\int\limits_{n}^{\infty}\psi(t)\frac{1}{\alpha(\psi;t)}dt.
\end{equation}
З рівності (\ref{eqq1}) випливає
\begin{equation}\label{eqq3}
\psi(n)n =
\int\limits_{n}^{\infty}\psi(t)\frac{1}{\alpha(\psi;t)}dt-\int\limits_{n}^{\infty}\psi(t)dt.
\end{equation}
Покажемо, що
\begin{equation}\label{eqq4}
\frac{1}{\alpha(\psi;t)}=1+\frac{1}{\alpha(g;t)}, \ t\geq1.
\end{equation}
Дійсно, оскільки $\psi(t)=g(t)t^{-1}$, то
$$
\frac{1}{\alpha(\psi;t)}=\frac{t|\psi'(t)|}{\psi(t)}=
\frac{t^{-1}g(t)+|g'(t)|}{g(t)t^{-1}}=
$$
$$
=1+\frac{t|g'(t)|}{g(t)}=1+\frac{1}{\alpha(g;t)}.
$$

Підставивши (\ref{eqq4}) в (\ref{eqq3}), отримаємо рівність
\begin{equation}\label{eqq5}
\psi(n)n=\int\limits_{n}^{\infty}\psi(t)\Big(1+\frac{1}{\alpha(g;t)}\Big)dt-\int\limits_{n}^{\infty}\psi(t)dt=
\int\limits_{n}^{\infty}\psi(t)\frac{1}{\alpha(g;t)}dt.
\end{equation}

З (\ref{eqq2}) i (\ref{eqq5})  випливає співвідношення
\begin{equation}\label{qq}
  \psi(n)n \leq
\frac{1}{\underline{\alpha}_{n}(g)}\int\limits_{n}^{\infty}\psi(t)dt \leq
\frac{1}{\underline{\alpha}_{n}(g)}
\sum\limits_{k=n}^{\infty}\psi(k).
\end{equation}
Нерівність (\ref{lemma_3}) доведено.

Нехай $g\in\mathfrak{M}_{C}$.
Оскільки $\mathfrak{M}_{C}\subset\mathfrak{M}_{0}$, то справедливість другої нерівності в (\ref{lemma_3_1}) випливає з (\ref{lemma_3}).

Врахувавши  (\ref{qq}), маємо
\begin{equation}\label{j2}
 \sum\limits_{k=n}^{\infty}\psi(k)\leq\psi(n)+\int\limits_{n}^{\infty}\psi(t)dt
 \leq\Big(\frac{1}{\underline{\alpha}_{n}(g)}\cdot\frac{1}{n}+1\Big)\int\limits_{n}^{\infty}\psi(t)dt.
\end{equation}

Тоді на підставі формул  (\ref{eqq5})--(\ref{j2}) одержуємо
$$
\psi(n)n \geq
\frac{1}{\overline{\alpha}_{n}(g)}\int\limits_{n}^{\infty}\psi(t)dt \geq
\frac{1}{\overline{\alpha}_{n}(g)}\Big(\frac{1}{\underline{\alpha}_{n}(g)}\cdot\frac{1}{n}+1\Big)^{-1}
\sum\limits_{k=n}^{\infty}\psi(k).
$$

Отже, співвідношення (\ref{lemma_3_1}), а отже і лему 3, доведено.

З (\ref{eq18}) і (\ref{lemma_3}) випливає нерівність
\begin{equation}\label{eq20}
  2\sum\limits_{k=n}^{2n-3}\psi(k)\Big(1-\frac{k}{2n-2}\Big)<\frac{1}{2\underline{\alpha}_{n}(g)}\sum\limits_{k=n}^{\infty}\psi(k),
\end{equation}

Об'єднуючи (\ref{for401}), (\ref{eq16}) і (\ref{eq20}),  отримаємо, що для довільних $l\in\mathbb{N}$ i $n\in\mathbb{N}$
\begin{equation}\label{w2}
|I_{2}|>\frac{1}{4}\Big|\cos\frac{\beta\pi}{2}\Big|
\Big(1-\frac{1}{2l}-\frac{1}{2\underline{\alpha}_{n}(g)}\Big)\sum\limits_{k=n}^{\infty}\psi(k).
\end{equation}

З іншої сторони, використовуючи
  твердження Д.1.1 з  \cite[с. 391]{Korn} та формулу (\ref{q3}), переконуємось, що для довільного $t_{n-1}(t)\in\mathcal{T}_{2n-1}$
  $$
  |I_{2}|\leq\|f_{D(l;n)}(t)-t_{n-1}(t)\|_{\infty}\|V_{D(l;n)}(t)-V_{n-1}(t)\|_{1}\leq
  $$
\begin{equation}\label{eq73}
\leq6\pi\|f_{D(l;n)}(t)-t_{n-1}(t)\|_{\infty}.
\end{equation}

З (\ref{w2}),  (\ref{eq73}) та умови  $\cos\frac{\beta\pi}{2}\neq0$ отримаємо
$$
E_{n}(f_{D(l;n)})_{\infty}=\inf\limits_{t_{n-1}\in \mathcal{T}_{2n-1}}\|f_{D(l;n)}(t)-t_{n-1}(t)\|_{\infty}\geq\frac{1}{6\pi}|I_{2}|\geq
$$
$$
 \geq\frac{1}{24\pi}\Big|\cos\frac{\beta\pi}{2}\Big|
\Big(1-\frac{1}{2n}-\frac{1}{2\underline{\alpha}_{n}(g)}\Big)\sum\limits_{k=n}^{\infty}\psi(k)\geq
$$
\begin{equation}\label{eq24}
  \geq\frac{1}{48\pi}\Big|\cos\frac{\beta\pi}{2}\Big|
\Big(1-\frac{1}{\underline{\alpha}_{1}(g)}\Big)\sum\limits_{k=n}^{\infty}\psi(k).
\end{equation}
Oб'єднуючи  (\ref{j3}) і (\ref{eq24}) отримуємо (\ref{theorem_3}). Теорему 3 доведено.

{\bf Теорема 4.} {\it Нехай ${\beta=2k+1},k\in\mathbb{Z}$, $\sum\limits_{k=1}^{\infty}\psi(k)<\infty$, ${\psi(t)=g(t)t^{-1}}$, $g\in\mathfrak{M}_{0}$, i
$$
\underline{\alpha}_{1}(g)=\inf\limits_{t\geq1}\alpha(g;t)>1.
$$
    Тоді для довільних $n\in \mathbb{N}$ мают місце нерівності
     \begin{equation}\label{theorem_4}
\frac{1}{360\pi}\Big(1-\frac{1}{\underline{\alpha}_{1}(g)}\Big)\psi(n)n \leq{ E}_{n}(L^{\psi}_{\beta,1})_{\infty} \leq{\cal E}_{n}(L^{\psi}_{\beta,1})_{\infty}\leq\Big(1+\frac{2}{\pi }\Big)\psi(n)n.
\end{equation}
 }
  {\bf Доведення теореми 4.}
 Спочатку знайдемо оцінку зверху для величини ${\cal E}_{n}(L^{\psi}_{\beta,1})_{\infty}$.
 З формули (\ref{for230}) при $\beta=2k+1, k\in\mathbb{Z}$, випливає
    \begin{equation}\label{for41}
{\cal E}_{n}(L^{\psi}_{\beta,1})_{\infty}\leq
\frac{1}{\pi}\big\|\Psi_{\beta,n}(t)\big\|_{\infty}=\frac{1}{\pi}\Big\|\sum\limits_{k=n}^{\infty}\psi(k)\sin kt\Big\|_{\infty}.
\end{equation}

Відомо
(див., наприклад, \cite[с. 611]{fiht})
\begin{equation}\label{for104}
\Big|\sum\limits_{j=1}^{k}\frac{\sin jx}{j}\Big|\leq \frac{\pi}{2}+1,  \ k\in \mathbb{N}, \ \ 0<x<2\pi.
\end{equation}
Застосувавши перетворення Абеля до суми $\sum\limits_{k=n}^{\infty}\psi(k)\sin kt$,
та використавши те, що $g(t)=\psi(t)t$  монотонно спадає, а також формулу (\ref{for104}), отримуємо
$$
\Big\|\sum\limits_{k=n}^{\infty}\psi(k)\sin kt\Big\|_{\infty}=
$$
$$
 =\Big\|\sum\limits_{k=n}^{\infty}(k\psi(k)-(k+1)\psi(k+1))\sum\limits_{j=1}^{k}
\frac{\sin jt}{j}-
n\psi(n)\sum\limits_{j=1}^{n-1}
\frac{\sin jt}{j}\Big\|_{\infty}\leq
$$
$$
\leq \Big(\frac{\pi}{2}+1\Big)\sum\limits_{k=n}^{\infty}|g(k)-g(k+1)|+\Big(\frac{\pi}{2}+1\Big)g(n)=
$$
\begin{equation}\label{for42}
=(\pi+2)g(n)=(\pi+2)\psi(n)n.
 \end{equation}

З формул (\ref{for41}) і (\ref{for42}) випливає оцінка зверху для величини ${\cal E}_{n}(L^{\psi}_{\beta,1})_{\infty}$ у випадку коли $\cos\frac{\beta\pi}{2}=0$.
А саме
\begin{equation}\label{for43}
 {\cal E}_{n}(L^{\psi}_{\beta,1})_{\infty}\leq \Big(1+\frac{2}{\pi }\Big)\psi(n)n.
\end{equation}

  Знайдемо оцінку знизу  величин  ${E}_{n}(L^{\psi}_{\beta,1})_{\infty}$.
Розглянемо функцію
\begin{equation}\label{q8}
  f^{*}_{n}(t)=\frac{1}{\pi}\int\limits_{-\pi}^{\pi}\Psi_{\beta}(\tau)\varphi^{*}_{n}(t-\tau)d\tau,
\end{equation}
де
\begin{equation}\label{for44}
\varphi^{*}_{n}(t)=\frac{-1}{5\pi n}\Big(\sum\limits_{k=1}^{n}k\sin kt+
\sum\limits_{k=n+1}^{2n}(2n+1-k)\sin kt\Big).
\end{equation}

Покажемо, що $\Vert\varphi^{*}_{n}\Vert_{1}\leq1$. Очевидно,
$$
\Big|\sum\limits_{k=1}^{n}k\sin kt+
\sum\limits_{k=n+1}^{2n}(2n+1-k)\sin kt\Big|\leq
$$
\begin{equation}\label{for45}
\leq \sum\limits_{k=1}^{n}k+\sum\limits_{k=n+1}^{2n}(2n+1-k)=n(n+1), \ 0\leq|t|\leq\pi.
\end{equation}
Застосовуючи перетворення Абеля до кожної з сум в (\ref{for44}), отримуємо
$$
\sum\limits_{k=1}^{n}k\sin kt+
\sum\limits_{k=n+1}^{2n}(2n+1-k)\sin kt=
$$
$$
=-\sum\limits_{k=1}^{n-1}\sum\limits_{j=1}^{k}\sin jt+n\sum\limits_{j=1}^{n}\sin jt+
\sum\limits_{k=n+1}^{2n-1}\sum\limits_{j=1}^{k}\sin jt+\sum\limits_{j=1}^{2n}\sin jt-n\sum\limits_{j=1}^{n}\sin jt=
$$
\begin{equation}\label{for47}
=-\sum\limits_{k=1}^{n-1}\sum\limits_{j=1}^{k}\sin jt+
\sum\limits_{k=n+1}^{2n}\sum\limits_{j=1}^{k}\sin jt.
\end{equation}

Використавши рівність
$$
\sum\limits_{k=0}^{N}\cos(kt+\gamma)=\cos\Big(\frac{N}{2}t +\gamma\Big)\sin
\frac{(N+1)t}{2}\cosec \frac{t}{2}, \ \ \gamma\in\mathbb{R}, \ 0<|t|\leq\pi,
$$
(див., наприклад, \cite[с. 43]{Gradshteyn})   з (\ref{for47})
отримуємо
$$
\Big|\sum\limits_{k=1}^{n}k\sin kt+
\sum\limits_{k=n+1}^{2n}(2n+1-k)\sin kt\Big|=
$$
$$
=\Big|-\sum\limits_{k=1}^{n-1}\frac{\cos\frac{t}{2}-\cos\left(k+\frac{1}{2}\right)t}{2\sin\frac{t}{2}}+
\sum\limits_{k=n+1}^{2n}\frac{\cos\frac{t}{2}-\cos\left(k+\frac{1}{2}\right)t}{2\sin\frac{t}{2}}\Big|=
$$
$$
=\frac{1}{2|\sin\frac{t}{2}|}
\Big|\sum\limits_{k=0}^{n-1}\cos\Big(k+\frac{1}{2}\Big)t-
\sum\limits_{k=n+1}^{2n}\cos\Big(k+\frac{1}{2}\Big)t\Big|=
$$
$$
=\frac{1}{2|\sin\frac{t}{2}|}
\Big|\sum\limits_{k=0}^{n-1}\cos\Big(k+\frac{1}{2}\Big)t-
\sum\limits_{k=0}^{2n}\cos\Big(k+\frac{1}{2}\Big)t+\sum\limits_{k=0}^{n}\cos\Big(k+\frac{1}{2}\Big)t\Big|\leq
$$
\begin{equation}\label{for48}
\leq\frac{3}{2(\sin\frac{t}{2})^{2}}, \ 0<|t|\leq \pi.
\end{equation}
Зі співвідношення (\ref{for48}) та нерівності
$$
\sin \frac{t}{2}\geq\frac{t}{\pi}, \ \
 0\leq t\leq\pi,
$$
маємо
\begin{equation}\label{for49}
\Big|\sum\limits_{k=1}^{n}k\sin kt+
\sum\limits_{k=n+1}^{2n}(2n+1-k)\sin kt\Big|\leq\frac{3\pi^2}{2t^{2}}, \ 0<|t|\leq \pi.
\end{equation}
З (\ref{for45}) i (\ref{for49}) випливає, що
$$
\left\Vert\varphi^{*}_{n}\right\Vert_{1}\leq \frac{1}{5\pi n}
\bigg( \int\limits_{|t|\leq\frac{\pi}{n}}n(n+1)dt+
\frac{3\pi^2}{2}\int\limits_{\frac{\pi}{n}\leq|t|\leq\pi}\frac{dt}{t^{2}}\bigg)=
$$
$$
=\frac{1}{5\pi n}
\big( 2\pi(n+1)+3\pi(n-1)\big)<1.
$$

Оскільки  $\Vert\varphi^{*}_{n}\Vert_{1}\leq1$ і  $\varphi^{*}_{n}\perp1$, то
$f^{*}_{n}\in L^{\psi}_{\beta,1}$.

З урахуванням формул (\ref{q8}), (\ref{for44}) та твердження (3.7.1) роботи \cite[с. 134]{Stepanets1} неважко переконатись, що для функції
$f^{*}_{n}$ вигляду (\ref{q8}) має місце рівність
\begin{equation}\label{eq26}
f^{*}_{n}(t)=\frac{1}{5\pi n}\Big(\sum\limits_{k=1}^{n}k\psi(k)\cos kt+
\sum\limits_{k=n+1}^{2n}(2n+1-k)\psi(k)\cos kt\Big).
\end{equation}

Розглянемо інтеграл
\begin{equation}\label{eq25}
  I_{3}=\int\limits_{-\pi}^{\pi}(f^{*}_{n}(t)-t_{n-1}(t))(V_{2n}(t)-V_{n}(t))dt,
\end{equation}
де $t_{n-1}\in\mathcal{T}_{2n-1}$, а $V_{m}(t)$ --- суми Валле Пуссена вигляду (\ref{for32}).
Використавши   твердження Д.1.1 з  \cite[с. 391]{Korn} та нерівність (\ref{q3}), отримаємо
$$
  I_{3}\leq\|f^{*}_{n}(t)-t_{n-1}(t)\|_{\infty}\|V_{2n}(t)-V_{n}(t)\|_{1}\leq
$$
\begin{equation}\label{eq27}
\leq6\pi\|f^{*}_{n}(t)-t_{n-1}(t)\|_{\infty}.
\end{equation}
Оскільки, згідно з (\ref{for32})
$$
 V_{2n}(t)-V_{n}(t)=
$$
$$
=\sum\limits_{k=n+1}^{2n}\cos kt+2\sum\limits_{k=2n+1}^{4n-1}\Big(1-\frac{k}{4n}\Big)\cos
kt-2\sum\limits_{k=n+1}^{2n-1}\Big(1-\frac{k}{2n}\Big)\cos
kt,
$$
то, враховуючи формули  (\ref{int_riv}),  (\ref{eq26}), (\ref{eq25}), та виконуючи елементарні перетворення, запишемо оцінку знизу
$$
 I_{3}=\int\limits_{-\pi}^{\pi}f^{*}_{n}(t)(V_{2n}(t)-V_{n}(t))dt=
$$
$$
=\frac{1}{5\pi n}\int\limits_{-\pi}^{\pi}\Big(\sum\limits_{k=1}^{n}k\psi(k)\cos kt+
\sum\limits_{k=n+1}^{2n}(2n+1-k)\psi(k)\cos kt\Big)\times
$$
$$
\times\Big(\sum\limits_{k=n+1}^{2n}\cos kt+2\sum\limits_{k=2n+1}^{4n-1}\Big(1-\frac{k}{4n}\Big)\cos kt-
2\sum\limits_{k=n+1}^{2n-1}\Big(1-\frac{k}{2n}\Big)\cos kt\Big)dt=
$$
$$
=\frac{1}{5 n}\Big(\sum\limits_{k=n+1}^{2n}(2n+1-k)\psi(k)-
2\sum\limits_{k=n+1}^{2n-1}(2n+1-k)\psi(k)\Big(1-\frac{k}{2n}\Big)\Big)=
$$
$$
=\frac{1}{5 n^{2}}\sum\limits_{k=n+1}^{2n}\psi(k)(2n+1-k)(k-n)\geq
$$
$$
=\frac{\psi(2n)}{5 n^{2}}\Big((2n+1)\sum\limits_{k=n+1}^{2n}(k-n)-\sum\limits_{k=n+1}^{2n}k^{2}+n\sum\limits_{k=n+1}^{2n}k\Big)=
$$
$$
=\frac{\psi(2n)}{5 n^{2}}\Big(\frac{n(n+1)(2n+1)}{2}-\frac{2n(2n+1)(4n+1)}{6}+
$$
$$
+\frac{n(n+1)(2n+1)}{6}+
\frac{n^{2}(3n+1)}{2}\Big)=
$$
$$
=\frac{\psi(2n)}{10 n}\Big(\frac{1}{3}n^{2}+n+\frac{2}{3}\Big)>\frac{1}{30}\psi(2n)n=
$$
\begin{equation}\label{eq28}
=\frac{1}{30}\psi(n)n\frac{\psi(2n)}{\psi(n)}=
  \frac{1}{60}\psi(n)n\frac{g(2n)}{g(n)}.
\end{equation}

Оскільки
$$
\frac{g(2n)}{g(n)}=\frac{g(2n)-g(n)}{ g(n)}+1>\frac{g'(n)n}{ g(n)}+1=1-\frac{1}{\alpha(g;n)}
\geq 1-\frac{1}{\underline{\alpha}_{1}(g)},
$$
то використовуючи  співвідношення
(\ref{eq27}), (\ref{eq28}), отримаємо, що для довільного $t_{n-1}\in\mathcal{T}_{2n-1}$ справедлива оцінка
$$
\|f^{*}_{n}(t)-t_{n-1}(t)\|_{\infty}\geq\frac{1}{6\pi}I_{3}\geq
$$
$$
\geq\frac{1}{6\pi} \frac{1}{60}\psi(n)n\frac{g(2n)}{g(n)}\geq
$$
$$
\geq
\frac{1}{360\pi}\Big(1-\frac{1}{\underline{\alpha}_{1}(g)}\Big)\psi(n)n.
$$
Тому, в силу довільності полінома $t_{n-1}\in\mathcal{T}_{2n-1}$, одержуємо оцінку
\begin{equation}\label{for105}
{ E}_{n}(L^{\psi}_{\beta,1})_{\infty}\geq \inf\limits_{t_{n-1}\in\mathcal{T}_{2n-1}}\|f^{*}_{n}(t)-t_{n-1}(t)\|_{\infty}
\geq
\frac{1}{360\pi}\Big(1-\frac{1}{\underline{\alpha}_{1}(g)}\Big)\psi(n)n.
\end{equation}
Об'єднуючи (\ref{for43}) і (\ref{for105}), отримуємо  (\ref{theorem_4}).
Теорему 4 доведено.

Прикладами функцій $\psi$, які задовольняють умови теорем 3--4  є функції:

1)  $\psi(t)=t^{-r}$, ${1<r<2}$;

2) $\psi(t)={t^{-1}\ln^{-\gamma}(t+K_{1})}$, $K_{1}\geq e^{\gamma}$, ${\gamma>1}$;

3) $\psi(t)={t^{-1}\ln^{-\gamma} (t+K_{1})
(\ln\ln (t+K_{2}))^{-\delta}}$, $\gamma\geq1$, $\delta>1$, ${K_{2}\geq K_{1}\geq e^{\gamma+\delta}}$,

 \noindent та інші.

Із теорем 3--4 безпосередньо випливає  твердження.

{\bf Наслідок 4. }{\it   Нехай $\sum\limits_{k=1}^{\infty}\psi(k)<\infty$, $\psi(t)=g(t)t^{-1}$, $g\in\mathfrak{M}_{0}$ і
$ \underline{\alpha}_{1}(g)>1$.
  Тоді,
  якщо $\beta\in\mathbb{R}$ такі, шо $\cos\frac{\beta\pi}{2}\neq0$, то
   для довільного $n\in \mathbb{N}$  справедливі порядкові оцінки
   \begin{equation}\label{ss2}
{ E}_{n}(L^{\psi}_{\beta,1})_{\infty}\asymp{\cal E}_{n}(L^{\psi}_{\beta,1})_{\infty}\asymp\sum\limits_{k=n}^{\infty}\psi(k);
\end{equation}
а якщо ж $\cos\frac{\beta\pi}{2}=0$, то
    \begin{equation}\label{ss3}
{ E}_{n}(L^{\psi}_{\beta,1})_{\infty}\asymp{\cal E}_{n}(L^{\psi}_{\beta,1})_{\infty}\asymp\psi(n)n.
\end{equation}
 }

 Якщо в умовах наслідку 4 $g\in\mathfrak{M}_{C}$, то
    згідно зі співвідношенням (\ref{lemma_3_1}),  має місце порядкова рівність
$$\sum\limits_{k=n}^{\infty}\psi(k)\asymp\psi(n)n.$$

Отже,  в цьому випадку має місце твердження.

{\bf Наслідок 5. }{\it  Нехай $\sum\limits_{k=1}^{\infty}\psi(k)<\infty$, $\psi(t)=g(t)t^{-1}$, $g\in\mathfrak{M}_{C}$,
$ \underline{\alpha}_{1}(g)>1$.
  Тоді   для довільних  $\beta\in\mathbb{R}$ справедливі порядкові оцінки }
  \begin{equation}\label{s1}
{ E}_{n}(L^{\psi}_{\beta,1})_{\infty}\asymp{\cal E}_{n}(L^{\psi}_{\beta,1})_{\infty}\asymp\psi(n)n.
\end{equation}

Справедливість порядкових оцінок (\ref{s1})  було встановлено раніше в \cite{Serdyuk_grabova}.

Наведемо наслідок з теорем 3--4 для функцій ${\psi(t)=t^{-1}\ln^{-\gamma}(t+e^{\gamma})}$, $\gamma>1$.

{\bf Наслідок 6. }{\it Нехай $\psi(t)=t^{-1}\ln^{-\gamma}(t+e^{\gamma}), \ \gamma>1$,  $\beta\in\mathbb{R}$,
 i $n\in\mathbb{N}, \ n\geq2$. Тоді}
 \begin{equation}\label{cons6}
  {E}_{n}(L^{\psi}_{\beta,1})_{\infty}\asymp{\cal E}_{n}(L^{\psi}_{\beta,1})_{\infty}\asymp {\left\{\begin{array}{cc}
\psi(n)n\ln n, \ \ & \cos\frac{\beta\pi}{2}\neq0, \\
\psi(n)n,  & \cos\frac{\beta\pi}{2}=0. \
  \end{array} \right.}
 \end{equation}
{\bf Доведення наслідку 6. } Покажемо, що для функцій ${\psi(t)=t^{-1}\ln^{-\gamma}(t+e^{\gamma })}$, $\gamma>1$, виконуються умови теорем 3--4.
  Дійсно, для них
 $$
 \sum\limits_{k=1}^{\infty}\psi(k)= \sum\limits_{k=1}^{\infty}\frac{1}{t\ln^{\gamma}(t+e^{\gamma})}<\infty
 $$
і
$$
\alpha(g;t)=\frac{\ln(t+e^{\gamma })}{\gamma}\frac{t+e^{\gamma}}{t}
>\frac{\ln(t+e^{\gamma })}{\gamma}>1.
$$
Якщо $\cos\frac{\beta\pi}{2}=0$, то порядкова оцінка (\ref{cons6}) безпосередньо випливає з (\ref{ss3}). Покажемо справедливість (\ref{cons6}) у випадку коли
 $\cos\frac{\beta\pi}{2}\neq0$.
Враховуючи нерівність (\ref{j2}) та монотонність функції $\psi(t)$ маємо
\begin{equation}\label{for1110}
\int\limits_{n}^{\infty}\psi(t)dt
\leq\sum\limits_{k=n}^{\infty}\psi(k)\leq\Big(\frac{1}{\underline{\alpha}_{n}(g)}\cdot\frac{1}{n}+1\Big)\int\limits_{n}^{\infty}\psi(t)dt
\leq2\int\limits_{n}^{\infty}\psi(t)dt.
\end{equation}

 Тоді з (\ref{ss2}) i (\ref{for1110}) одержуємо
$$
{ E}_{n}(L^{\psi}_{\beta,1})_{\infty}\asymp {\cal E}_{n}(L^{\psi}_{\beta,1})_{\infty}\asymp
$$
$$
 \asymp\sum\limits_{k=n}^{\infty}\psi(k)\asymp\int\limits_{n}^{\infty}\psi(t)dt\asymp
 \int\limits_{n}^{\infty}\frac{dt}{t\ln^{\gamma}(t+e^{\gamma})}\asymp
$$
$$
\asymp\ln^{1-\gamma}n
\asymp\psi(n)n\ln n, \ n\geq2.
$$
Наслідок 6  доведено.

E-mail: \href{mailto:tania_stepaniuk@ukr.net}{tania$_{-}$stepaniuk@ukr.net}


\newpage
\begin{thebibliography}{10}


\bibitem{Stepanets1}{\sc Степанец А.И.}
 \emph{Методы теории приближений}: В 2 ч. // Праці Інституту математики НАН України --- Киев: Ин-т
математики НАН Украины, 2002. --- {\bf 40}. --- Ч.І. --- 427 с.

\bibitem{Zigmund2}{\sc Зигмунд А.}
 \emph{Тригонометрические ряды.}  В 2 т.~--- М.: Мир,
1965.~--- Т. ІІ.~ --- 538~с.

\bibitem{Korn}{\sc Корнейчук Н.П.}
 \emph{Точные константы в теории приближения.}  --- М.:
Наука, 1987. ---  424~с.

\bibitem{T}{\sc Temlyakov V.N.}
 \emph{Approximation of Periodic Function}: NY: Nova
Science Publichers, Inc. --- 1993. --- 419p.

 \bibitem{Serdyuk_grabova}{\sc  Сердюк А.С., Грабова У.З.}
 Порядкові оцінки найкращих
наближень і наближень сумами Фур'є  класів $(\psi,\beta)$ --
диференційовних функцій
// Укр. мат. журн. --- 2013.
--- {\bf 65}, №9. --- С. 1186 -- 1197.

\bibitem{Step monog 1987}{\sc Степанец А.И.}
 \emph{Классификация и приближение периодических функций.} --- Киев: Наук. думка~--- 1987.~--- 268 c.


\bibitem{S_S}{\sc  Сердюк А.С., Степанюк Т.А.}
Порядкові оцінки найкращих
наближень і наближень сумами Фур'є класів нескінченно
диференційовних функцій// Теорія наближення функцій та суміжні питання: Зб. праць Ін--ту математики НАН України.
2013. --- {\bf 10}, №1.
 --- С. 255--282.

 \bibitem{Serdyuk2013}{\sc  Сердюк А.С., Соколенко І.В.}
Наближення лінійними методами
класів  $(\psi, \overline{\beta})$--диференційовних функцій
// Теорія наближення функцій та суміжні питання: Зб. праць Ін-ту математики НАН України. --- Київ: Інстиут математики НАН України, 2013. --- {\bf 10}, №1.
 --- С. 245--254.

 \bibitem{serdyuk2004zbirnyk}{\sc  Сердюк А.С.}
 Про один лінійний метод наближення періодичних
функцій// Проблеми теорії наближеня функцій та суміжні питанння: Зб. праць Ін-ту матем. НАН України.
--- Київ: Інстиут математики НАН України, --- 2004. --- Т.1, №1.
 --- С. 294--336.

  \bibitem{Serdyuk2002}{\sc  Сердюк А.С.}
  Про найкраще наближення на класах згорток
періодичних функцій
// Теорія наближення та її застосування: Пр.  Ін-ту математики НАН України. Т.41. ---
Київ: Ін-т математики НАН України, 2002.
 --- С.~168--189.

\bibitem{Serdyuk2005}{\sc  Сердюк А.С.}
  Найкращі наближення і поперечники класів згорток періодичних функцій високої гладкості
// Укр. мат. журн. --- 2005.
--- {\bf 57}, №7. --- С.946 -- 971.


 \bibitem{Serduyk_Stepaniuk}{\sc  Serdyuk A.S., Stepaniuk T.A.}
 Order estimates of the best approximations and approximations of
 Fourier sums of classes of convolutions of periodic functions of not high smoothness in uniform metric //
  Arxiv preprint, arXiv:1403.5311, 2014. --- 20~p.


\bibitem{Zigmund1}{\sc Зигмунд А.}
 \emph{Тригонометрические ряды.}  В 2 т.~--- М.: Мир,
1965.~--- Т. І.~ --- 615~с.


\bibitem{Gradshteyn}{\sc Градштейн И.С., Рыжик И.М.}
 \emph{Таблицы интегралов, сумм, рядов и
произведений.} --- М.: Физматиз, 1962. --- 1100 с.


\bibitem{fiht}{\sc Фихтенгольц Г.М.}
 \emph{Курс дифференциального и интегрального исчисления.}  В 3 т.~--- М.: Наука,
1969.~--- Т.ІІІ.~ --- 656 с.










\end{thebibliography}
\end{document}